\documentclass[12pt, reqno]{amsart}
\usepackage{amsmath, amsthm, amscd, amsfonts, amssymb, graphicx, color}
\usepackage[bookmarksnumbered, colorlinks, plainpages]{hyperref}

\newtheorem{theorem}{Theorem}[section]
\newtheorem{lemma}[theorem]{Lemma}
\newtheorem{proposition}[theorem]{Proposition}
\newtheorem{corollary}[theorem]{Corollary}
\theoremstyle{definition}
\newtheorem{definition}[theorem]{Definition}
\newtheorem{example}[theorem]{Example}

\theoremstyle{remark}
\newtheorem{remark}[theorem]{Remark}
\numberwithin{equation}{section}

\allowdisplaybreaks
\begin{document}

\title[Exact Continuous frames in Hilbert $C^*$-modules]
{Exact Continuous frames in Hilbert $C^*$-modules}

\author[H.Ghasemi]{Hadi Ghasemi }
\address{Hadi Ghasemi \\ Department of Mathematics and Computer
Sciences, Hakim Sabzevari University, Sabzevar, P.O. Box 397, IRAN}
\email{ \rm h.ghasemi@hsu.ac.ir}
\author[T.L. Shateri]{Tayebe Lal Shateri }
\address{Tayebe Lal Shateri \\ Department of Mathematics and Computer
Sciences, Hakim Sabzevari University, Sabzevar, P.O. Box 397, IRAN}
\email{ \rm  t.shateri@hsu.ac.ir; shateri@ualberta.ca}
\thanks{*The corresponding author:
t.shateri@hsu.ac.ir; shateri@ualberta.ca (Tayebe Lal Shateri)}
 \maketitle

\begin{abstract}
In the present paper, we investigate some properties of  continuous frames in Hilbert $C^*$-modules. In particular, we give the requirements so that by removing some elements of a continuous frame, it does not remain a continuous frame and when the remaining set still remains a continuous frame. Finally, we consider the requirements under which a continuous frame is not a Riesz-type frame.
\vskip 3mm
\end{abstract}

\section{Introduction And Preliminaries}
Duffin and Schaeffer \cite{DS} introduced the concept of a frame to study some deep problems in the nonharmonic Fourier series. After the fundamental paper by Daubechies et al. \cite{DG}, various generalizations of frames were developed. Now, frames are useful in some areas such as signal and image processing, filter bank theory, data compression and sampling theory, etc. We refer to \cite{CH} for an introduction to frame theory in Hilbert spaces and its applications. 
 
In 2000, Frank-Larson \cite{FL} introduced the notion of frames in Hilbert $C^*$-modules as a generalization of frames in Hilbert spaces. It is well known that Hilbert $C^*$-modules are generalizations of Hilbert spaces by allowing the inner product to take values in a $C^*$-algebra rather than in the field of complex numbers. The theory of Hilbert $C^*$-modules has applications in the study of locally compact quantum groups, complete maps between -algebras, non-commutative geometry, and KK-theory. There are some
differences between Hilbert $C^*$-modules and Hilbert spaces. For example, we know that the Riesz representation theorem for continuous linear functionals on Hilbert spaces does not
extend to Hilbert $C^*$-modules \cite{OLS}, also any closed subspace in a Hilbert space has an orthogonal complement, but it is not true for Hilbert $C^*$-module \cite{MA}. Moreover, we know that every bounded operator on a Hilbert space has an adjoint, while there are bounded operators on Hilbert $C^*$-modules which do not have any \cite{MT}. Thus it is more difficult to make a discussion of the theory of Hilbert $C^*$-modules than those for Hilbert spaces. This makes the
study of the frames for Hilbert $C^*$-modules important and interesting. We refer the readers to \cite{LAN}, for more details on Hilbert $C^*$-modules. 
The theory of frames has been extended from Hilbert spaces to Hilbert $C^*$-modules, see \cite{RR,SH,XZ}.

The paper is organized as follows. First, we recall the basic definitions and some notations about Hilbert $C^*$-modules, and we also give some properties of them which we will use in the later sections. Also, we recall the definition of continuous frames in Hilbert $C^*$-modules. In Section 2, we deal to the duals of continuous frames in Hilbert $C^*$-modules and prove some important properties of them. In particular, we show when can remove some elements from a continuous frame so that the remaining set is not a continuous frame and when the remaining set still remains a continuous frame. Finally, we consider the requirements under which a continuous frame has more duals.

First, we recall some definitions and basic properties of Hilbert $C^*$-modules. We give only a brief introduction to the theory of Hilbert $C^*$-modules to make our explanations self-contained. For comprehensive accounts, we refer to  \cite{LAN,OLS}.
We now give the definition of Hilbert $C^*$-modules. Throughout this paper, $\mathcal A$ shows a unital $C^*$-algebra.
%-------------------------------------------------------------------------------------------------------------------------------------------------------
\begin{definition}
A \textit{pre-Hilbert module} over unital $C^*$-algebra $\mathcal A$ is a complex vector space $U$ which is also a left $\mathcal A$-module equipped with an $\mathcal A$-valued inner product $\langle .,.\rangle :U\times U\to \mathcal A$ which is $\mathbb C$-linear and $\mathcal A$-linear in its first variable and satisfies the following conditions:\\
$(i)\; \langle f,f\rangle \geq 0$,\\
$(ii)\; \langle f,f\rangle =0$  iff $f=0$,\\
$(iii)\; \langle f,g\rangle ^*=\langle g,f\rangle ,$\\
$(iv)\; \langle af,g\rangle=a\langle f,g\rangle ,$\\
for all $f,g\in U$ and $a\in\mathcal A$.
\end{definition}
%---------------------------------------------------------------------------------------------------------------------------------------------------------
A pre-Hilbert $\mathcal A$-module $U$ is called \textit{Hilbert $\mathcal A$-module} if $U$ is complete with respect to the topology determined by the norm $\|f\|=\|\langle f,f\rangle \|^{\frac{1}{2}}$.

By \cite[Example 2.46]{JI}, if $\mathcal A$ is a $C^*$-algebra, then it is a Hilbert $\mathcal A$-module with respect to the inner product
$$\langle a,b\rangle =ab^*,\quad (a,b\in \mathcal A).$$
%---------------------------------------------------------------------------------------------------------------------------------------------------------
\begin{example}
\cite[Page 237]{OLS} Let $l^2(\mathcal A)$ be the set of all sequences $\{a_n\}_{n\in \mathbb N}$ of elements of a $C^*$-algebra $\mathcal A$ such that the series $\sum_{n=1}^{\infty}a_na_n^*$ is convergent in $\mathcal A$. Then $l^2(\mathcal A)$ is a Hilbert $\mathcal A$-module with respect to the pointwise operations and inner product defined by
\begin{equation*}
\langle \{a_n\}_{n\in \mathbb N},\{b_n\}_{n\in \mathbb N}\rangle =\sum_{n=1}^{\infty}a_nb_n^*.
\end{equation*}
\end{example}
%---------------------------------------------------------------------------------------------------------------------------------------------------------
In the following lemma the \textit{Cauchy-Schwartz inequality} reconstructed in Hilbert $C^*$-modules.
%---------------------------------------------------------------------------------------------------------------------------------------------------------
\begin{lemma}
\cite[Lemma 15.1.3]{OLS} (\textbf{Cauchy-Schwartz inequality}) Let $U$ be a Hilbert $C^*$-modules over a unital $C^*$-algebra $\mathcal A$. Then
\begin{equation*}
\|\langle f,g\rangle \|^{2}\leq\|\langle f,f\rangle \|\;\|\langle g,g\rangle \|
\end{equation*}
for all $f,g\in U$.
\end{lemma}
%---------------------------------------------------------------------------------------------------------------------------------------------------------
\begin{proposition}\label{NO}
\cite[Proposition 1.2.4]{MT} Let $U$ be a Hilbert $C^*$-modules over a unital $C^*$-algebra $\mathcal A$. Then\\
(i)\;$\| af\|\leq\| a\|\| f\|$\\
(ii)\;$\langle f,g\rangle\langle g,f\rangle=\| g\|\langle f,f\rangle$\\
for all $f,g\in U$ and $a\in\mathcal A$.
\end{proposition}
%-------------------------------------------------------------------------------------------------------------------------------------------
In the following, we assume that $\mathcal A$ is a unital $C^*$-algebra and $U$ is a Hilbert $C^*$-module over $\mathcal A$ and $(\Omega ,\mu)$ is a measure space. We give the definition of continuous frames for Hilbert $\mathcal A$-modules and some properties of them from \cite{GHSH}.
%-------------------------------------------------------------------------------------------------------------------------------------------

Let $(\Omega ,\mu)$ be a measure space and $\mathcal A$ is a unital $C^*$-algebra. Consider
\begin{equation*}
L^{2}(\Omega ,A)=\lbrace\varphi :\Omega \to A\quad ;\quad \Vert\int_{\Omega}\vert(\varphi(\omega))^{*}\vert^{2} d\mu(\omega)\Vert<\infty\rbrace.
\end{equation*}
For any $\varphi ,\psi \in L^{2}(\Omega , A)$, the inner product is defined by $\langle \varphi ,\psi\rangle = \int_{\Omega}\langle\varphi(\omega),\psi(\omega)\rangle d\mu(\omega)$ and the norm is defined by $\|\varphi\|=\|\langle \varphi,\varphi\rangle \|^{\frac{1}{2}}$. It was shown in \cite{LAN} $L^{2}(\Omega , A)$ is a Hilbert $\mathcal A$-module.
%-------------------------------------------------------------------------------------------------------------------------------------------
\begin{definition}
A mapping $F:\Omega \to U$ is called a continuous frame for $U$ if\\
$(i)\; F$ is weakly-measurable, i.e, for any $f\in U$, the mapping $\omega\longmapsto\langle f,F(\omega)\rangle $ is measurable on $\Omega$.\\
$(ii)$ There exist constants $A,B>0$ such that
\begin{equation}\label{eq1}
A\langle f,f\rangle \leq \int_{\Omega}\langle f,F(\omega)\rangle \langle F(\omega),f\rangle d\mu(\omega)\leq B\langle f,f\rangle  ,\quad (f\in U). 
\end{equation}
\end{definition}
%-------------------------------------------------------------------------------------------------------------------------------------------
The constants $A,B$ are called \textit{lower} and \textit{upper} frame bounds, respectively. The mapping $F$ is called \textit{Bessel} if the right inequality in \eqref{eq1} holds and is called \textit{tight} if  $A=B$.\\
%-------------------------------------------------------------------------------------------------------------------------------------------
A continuous frame $F:\Omega \to U$ is called \textit{exact} if for every measurable subset $\Omega_{1}\subseteq\Omega$ with
$0<\mu(\Omega_{1})<\infty$, the mapping $F:\Omega\backslash\Omega_{1} \to U$ is not a continuous frame for $U$.\\
%-------------------------------------------------------------------------------------------------------------------------------------------
Now we give definitions of the pre-frame operator and the continuous frame operator for a continuous frame in Hilbert $C^{\ast}$-modules.
%-------------------------------------------------------------------------------------------------------------------------------------------

Let $F:\Omega \to U$ be a Bessel mapping . Then \\
(i)\;The \textit{synthesis operator} or \textit{pre-frame operator} $T_{F}:L^{2}(\Omega , A)\;\to U$ weakly defined by
\begin{equation}
\langle T_{F}\varphi ,f\rangle =\int_{\Omega}\varphi(\omega)\langle F(\omega),f\rangle d\mu(\omega),\quad (f\in U).
\end{equation}
(ii)\; The adjoit of $T$, called The \textit{analysis operator} $T^{\ast}_{F}:U\;\to L^{2}(\Omega , A)$ is defined by
\begin{equation}
(T^{\ast}_{F}f)(\omega)=\langle f ,F(\omega)\rangle\quad (\omega\in \Omega).
\end{equation}
%-------------------------------------------------------------------------------------------------------------------------------------------
The frame operator $S_{F}:U\;\to U$ for a continuous frame $F:\Omega \to U$ is weakly defined by
\begin{equation}
\langle S_{F}f ,f\rangle =\int_{\Omega}\langle f,F(\omega)\rangle\langle F(\omega),f\rangle d\mu(\omega),\quad (f\in U).
\end{equation}
%-------------------------------------------------------------------------------------------------------------------------------------------
In \cite{GHSH} it is shown that the pre-frame operator $T_{F}:L^{2}(\Omega , A)\;\to U$ is well defined, surjective, adjointable and  $\| T\|\leq\sqrt{B}$ . Moreover the analysis operator $T^{\ast}_{F}:U\;\to L^{2}(\Omega , A)$ is injective and has closed range. Also $S=T T^{\ast}$ is positive, adjointable, self-adjoint and invertible and $\| S\|\leq B$.
%-------------------------------------------------------------------------------------------------------------------------------------------
 \section{Main results}
In this section, we give the consept of duals of continuous frames in Hilbert $C^{\ast}$-modules and prove some important properties of continuous frames and their duals. In particular, we show when can remove some elements from a continuous frame so that the remaining set is not a continuous frame and when the remaining set still remains a continuous frame.
 %-------------------------------------------------------------------------------------------------------------------------------------------
\begin{definition}
\cite{GHSH} Let $F:\Omega \to U$ be a continuous Bessel mapping. A continuous Bessel mapping $G:\Omega \to U$ is called a \textit{dual} for $F$ if
\begin{equation*}
f= \int_{\Omega}\langle f,G(\omega)\rangle F(\omega)d\mu(\omega)\;\;\;\;\;(f\in U)
\end{equation*}
or
\begin{equation}\label{eq2}
\langle f,g\rangle= \int_{\Omega}\langle f,G(\omega)\rangle\langle F(\omega),g\rangle d\mu(\omega)\;\;\;\;\;(f,g\in U)
\end{equation}
In this case $(F,G)$ is called a \textit{dual pair}. If $T_{F}$ and $T_{G}$ denote the synthesis operators of $F$ and $G$, respectively, then \eqref{eq2} is equivalent to $T_{F}T^{*}_{G}=I_{U}$.
\end{definition}
%-------------------------------------------------------------------------------------------------------------------------------------------
 It is easy to see that $S^{-1}F$ is a dual for $F$, which is called \textit{canonical dual}.
%-------------------------------------------------------------------------------------------------------------------------------------------
\begin{definition}
\cite{GHSH} Let $F:\Omega \to U$ be a continuous frame for Hilbert $C^{\ast}$-module $U$. If $F$ has only one dual, we call $F$ a \textit{Riesz-type frame}.
\end{definition}
%-------------------------------------------------------------------------------------------------------------------------------------------
In \cite{GHSH}, it is shown that a continuous frame $F$ for Hilbert $C^{\ast}$-module $U$ over a unital $C^*$-algebra $\mathcal A$ is a Riesz-type frame if and only if the analysis operator $T^{*}_{F}:U\to L^{2}(\Omega ,A)$ is onto.
%-------------------------------------------------------------------------------------------------------------------------------------------
\begin{proposition}
Let $F:\Omega \to U$ be a Riesz-type frame for Hilbert $C^{\ast}$-module $U$. Then $F(\omega)\neq 0$ for every $\omega\in\Omega$.
\end{proposition}
\begin{proof}
Let $G$ be the canonical dual of $F$ and $\omega_{0}\in\Omega$ such that $F(\omega_{0})=0$ . Define $G_{1}:\Omega \to U$ where $G_{1}(\omega_{0})\neq 0$ and $G_{1}(\omega)=G(\omega)$ for all $\omega\neq\omega_{0}$. Then $G_{1}$ is a continuous Bessel mapping and
\begin{align*}
f & =\int_{\Omega}\langle f,G(\omega)\rangle F(\omega) d\mu(\omega)\\
& =\int_{\lbrace\omega_{0}\rbrace}\langle f,G(\omega)\rangle F(\omega) d\mu(\omega) +\int_{\Omega\setminus\lbrace\omega_{0}\rbrace}\langle f,G(\omega)\rangle F(\omega) d\mu(\omega)\\
& =\int_{\lbrace\omega_{0}\rbrace}\langle f,G_{1}(\omega)\rangle F(\omega) d\mu(\omega) +\int_{\Omega\setminus\lbrace\omega_{0}\rbrace}\langle f,G_{1}(\omega)\rangle F(\omega) d\mu(\omega)\\
& =\int_{\Omega}\langle f,G_{1}(\omega)\rangle F(\omega) d\mu(\omega).
\end{align*}
Hence $G_{1}$ is a dual of $F$ and $f$ is not Riesz-type.
\end{proof}
%-------------------------------------------------------------------------------------------------------------------------------------------
The following theorem shows that for each dual pair $(F,G)$ for Hilbert $C^{\ast}$-module $U$, both $F$ and $G$ are continuous frames for $U$.
%-------------------------------------------------------------------------------------------------------------------------------------------
\begin{theorem}
Let $F:\Omega \to U$ be a continuous Bessel mapping for Hilbert $C^{\ast}$-module $U$ with bound $B_{F}$. If there exists a continuous Bessel mapping $G:\Omega \to U$ with bound $B_{G}$ such that for any $f,g\in U$ ,
\begin{equation*}
\langle f,g\rangle=\int_{\Omega}\langle f,G(\omega)\rangle\langle F(\omega),g\rangle d\mu(\omega),
\end{equation*}
then $F$ is a continuous frame for $U$.
\end{theorem}
\begin{proof}
For any $f\in U$ we have,
\begin{align*}
\Vert\langle f,f\rangle\Vert^{2} & =\Vert\int_{\Omega}\langle f,G(\omega)\rangle\langle F(\omega),f\rangle d\mu(\omega)\Vert^{2}\\
& =\Vert\langle\lbrace\langle f,G(\omega)\rangle\rbrace_{\omega\in\Omega},\lbrace\langle f,F(\omega)\rangle\rbrace_{\omega\in\Omega}\rangle\Vert^{2}\\
& \leq\Vert\lbrace\langle f,G(\omega)\rangle\rbrace_{\omega\in\Omega}\Vert^{2}\Vert\lbrace\langle f,F(\omega)\rangle\rbrace_{\omega\in\Omega}\Vert^{2}\\
& =\Vert\int_{\Omega}\langle f,G(\omega)\rangle\langle G(\omega),f\rangle d\mu(\omega)\Vert\Vert\int_{\Omega}\langle f,F(\omega)\rangle\langle F(\omega),f\rangle d\mu(\omega)\Vert\\
& \leq B_{G}\Vert\langle f,f\rangle\Vert\Vert\int_{\Omega}\langle f,F(\omega)\rangle\langle F(\omega),f\rangle d\mu(\omega)\Vert.
\end{align*}
Hence
\begin{equation*}
B^{-1}_{G}\Vert\langle f,f\rangle\Vert\leq\Vert\int_{\Omega}\langle f,F(\omega)\rangle\langle F(\omega),f\rangle d\mu(\omega)\Vert\leq B_{F}\Vert\langle f,f\rangle\Vert.
\end{equation*}
This shows that $F$ is a continuous frame for $U$, by \cite[Theorem 2.13]{GHSH}
\end{proof}
%-------------------------------------------------------------------------------------------------------------------------------------------
\begin{remark}
By the existence of a canonical dual frame for each continuous frame, the converse of the above theorem is obvious.
\end{remark}
%-------------------------------------------------------------------------------------------------------------------------------------------
\begin{example} 
Let
$U=\Big\{
\begin{pmatrix}
a&0&0\\
0&0&b
\end{pmatrix}
: a,b\in \mathbb C\Big\}$, and $\mathcal A=\Big\{
\begin{pmatrix}
x&0\\
0&y
\end{pmatrix}
: x,y\in \mathbb C\Big\}$ which is an unital $C^*$-algebra. We define the inner product
\begin{equation*}
\begin{array}{ll}
\langle .,.\rangle:U\times U\;\to \quad \mathcal A \\
\qquad\; (M,N)\longmapsto M(\overline{N})^t.
\end{array}
\end{equation*}
This inner product makes $U$ a Hilbert $C^*$-module over $\mathcal A$. Suppose that $(\Omega ,\mu)$ is a measure space where  $\Omega=[0,1]$ and $\mu$ is the Lebesgue measure. Consider the continuous Bessel mappings $F,G:\Omega\to U$ defined by 
$F(\omega)=\begin{pmatrix}
2\omega&0&0\\
0&0&2\omega
\end{pmatrix}$ and $G(\omega)=\begin{pmatrix}
\dfrac{3}{2}\omega&0&0\\
0&0&\dfrac{3}{2}\omega
\end{pmatrix}$, for any $\omega\in \Omega$.\\ 
For each $f=
\begin{pmatrix}
a&0&0\\
0&0&b
\end{pmatrix}
\in U$, we have
\begin{align*}
\int_{\Omega}\langle f,F(\omega)\rangle G(\omega) d\mu(\omega) &=\int_{[0,1]}
\begin{pmatrix}
2\omega a&0\\
0&2\omega b
\end{pmatrix}\begin{pmatrix}
\dfrac{3}{2}\omega&0&0\\
0&0&\dfrac{3}{2}\omega
\end{pmatrix} d\mu(\omega)\\
&=\begin{pmatrix}
a&0&0\\
0&0&b
\end{pmatrix}\int_{[0,1]}3\omega^{2} d\mu(\omega)=\begin{pmatrix}
a&0&0\\
0&0&b
\end{pmatrix} =f.\\
\end{align*}
Therefore $G$ is a dual of $F$.
Moreover, both $F$ and $G$ are continuous tight frames for $U$ where $A_{F}=B_{F}=\dfrac{4}{3}$ is the frame bound of $F$ and $A_{G}=B_{G}=\dfrac{3}{4}$ is the frame bound of $G$.
\end{example}
%-------------------------------------------------------------------------------------------------------------------------------------------
\begin{example} 
Assume that $\mathcal A=\Big\{\begin{pmatrix}
a&0\\
0&b
\end{pmatrix}: x,y\in \mathbb C\Big\}$ which is an unital $C^*$-algebra. We define the inner product
\begin{equation*}
\begin{array}{ll}
\langle .,.\rangle:\mathcal A\times \mathcal A\;\to \quad \mathcal A \\
\qquad\; (M,N)\longmapsto M(\overline{N})^t.
\end{array}
\end{equation*}
$\mathcal A$ with this inner product is a Hilbert $C^*$-module over itself. Suppose that $(\Omega ,\mu)$ is a measure space where  $\Omega=[0,1]$ and $\mu$ is the Lebesgue measure. Consider the continuous Bessel mappings $F,G:\Omega\to \mathcal A$ defined by 
$F(\omega)=\begin{pmatrix}
2\omega&0\\
0&\omega -1
\end{pmatrix}$ and $G(\omega)=\begin{pmatrix}
\dfrac{3}{2}\omega&0\\
0&\omega-\dfrac{7}{3}
\end{pmatrix}$, for any $\omega\in \Omega$.\\ 
For each $f=
\begin{pmatrix}
a&0\\
0&b
\end{pmatrix}
\in \mathcal A$, we have
\begin{align*}
\int_{\Omega}\langle f,F(\omega)\rangle G(\omega) d\mu(\omega) &=\int_{[0,1]}
\begin{pmatrix}
2\omega a&0\\
0&(\omega-1) b
\end{pmatrix}\begin{pmatrix}
\dfrac{3}{2}\omega&0\\
0&\omega-\dfrac{7}{3}
\end{pmatrix} d\mu(\omega)\\
&=\begin{pmatrix}
a&0\\
0&b
\end{pmatrix}\int_{[0,1]}
\begin{pmatrix}
3\omega^{2} &0\\
0&\omega^{2}-\dfrac{10}{3}\omega+\dfrac{7}{3}
\end{pmatrix} d\mu(\omega)\\
&=\begin{pmatrix}
a&0\\
0&b
\end{pmatrix}\begin{pmatrix}
1&0\\
0&1
\end{pmatrix} =\begin{pmatrix}
a&0\\
0&b
\end{pmatrix}=f.\\
\end{align*}
Therefore $F$ is a dual of $G$. It is easy to see that both $F$ and $G$ are continuous frames for $\mathcal A$ where $\dfrac{1}{3},\dfrac{4}{3}$ are the frame bounds of $F$ and $\dfrac{3}{4},\dfrac{31}{9}$ are the frame bounds of $G$.
\end{example}
%-------------------------------------------------------------------------------------------------------------------------------------------
In \cite{JI}, it has stated an important property of a given frame for Hilbert $C^{\ast}$-module $U$ and showed when it can removed some elements from a frame so that the set still remains a frame. Now we generalize these properties to the situation of continuous frames for Hilbert $C^{\ast}$-module $U$ over a unital $C^*$-algebra $\mathcal A$.\\

We need the following lemma for the proof of the next theorem.              
%-------------------------------------------------------------------------------------------------------------------------------------------
\begin{lemma}\label{psi}
Let $F:\Omega \to U$ be a continuous frame for Hilbert $C^{\ast}$-module $U$ over a unital $C^*$-algebra $\mathcal A$ with frame operator $S$ and $\omega_{0}\in\Omega$. Then the mapping
\begin{align*}
\psi: & \Omega\longrightarrow\mathcal A\\
& \omega\longmapsto\langle F(\omega_{0}),S^{-1}F(\omega)\rangle
\end{align*}
belongs to $L^{2}(\Omega , \mathcal A)$ and
\begin{equation*}
 \psi(\omega_{0})=\int_{\Omega\setminus\lbrace\omega_{0}\rbrace}\psi(\omega)\psi(\omega)^{*} d\mu(\omega)+(\psi(\omega_{0}))^{2}\mu(\lbrace\omega_{0}\rbrace).
\end{equation*}
\end{lemma}

\begin{proof}
Assume that the upper frame bound of $F$ is $B_{F}$. Then
\begin{align*}
\int_{\Omega}\psi(\omega)\psi(\omega)^{*} d\mu(\omega) &= \int_{\Omega}\langle F(\omega_{0}),S^{-1}F(\omega)\rangle\langle S^{-1}F(\omega),F(\omega_{0}) d\mu(\omega)\\
& = \int_{\Omega}\langle S^{-1}F(\omega_{0}),F(\omega)\rangle\langle F(\omega),S^{-1}F(\omega_{0})\rangle d\mu(\omega)\\
& \leq B_{F}\langle S^{-1}F(\omega_{0}),S^{-1}F(\omega_{0})\rangle.
\end{align*}
Hence
\begin{equation*}
\Vert\int_{\Omega}\psi(\omega)\psi(\omega)^{*} d\mu(\omega)\Vert\leq B_{F}\Vert S^{-1}F(\omega_{0})\Vert^{2}<\infty,
\end{equation*}
i.e. $\psi\in L^{2}(\Omega ,\mathcal A)$. Also
\begin{align*}
\psi(\omega_{0}) &=\langle F(\omega_{0}),S^{-1}F(\omega_{0})\rangle\\
&= \langle\int_{\Omega}\langle F(\omega_{0}),S^{-1}F(\omega)\rangle F(\omega) d\mu(\omega),S^{-1}F(\omega_{0})\rangle\\
&= \int_{\Omega}\psi(\omega) \langle F(\omega) ,S^{-1}F(\omega_{0})\rangle d\mu(\omega)\\
&= \int_{\Omega\setminus\lbrace\omega_{0}\rbrace}\psi(\omega) \langle F(\omega) ,S^{-1}F(\omega_{0})\rangle d\mu(\omega)+\psi(\omega_{0})\langle F(\omega_{0}) ,S^{-1}F(\omega_{0})\rangle\mu(\lbrace\omega_{0}\rbrace)\\
&= \int_{\Omega\setminus\lbrace\omega_{0}\rbrace}\psi(\omega)\psi(\omega)^{*} d\mu(\omega)+(\psi(\omega_{0}))^{2}\mu(\lbrace\omega_{0}\rbrace).
\end{align*}
Moreover, $\psi(\omega_{0})$ is self-adjoint.
\end{proof}                
 %-------------------------------------------------------------------------------------------------------------------------------------------
The following theorem shows when does removing some elements from a continuous frame cause that the remaining set is not a continuous frame.              
 %-------------------------------------------------------------------------------------------------------------------------------------------
\begin{theorem}\label{c-f-less-no}
Let $F:\Omega \to U$ be a continuous frame for Hilbert $C^{\ast}$-module $U$ over a unital $C^*$-algebra $\mathcal A$ with frame operator $S$. Let $1_{\mathcal A}$ be the identity element of $\mathcal A$ and $\omega_{0}\in\Omega$ such that $1_{\mathcal A}-\langle F(\omega_{0}),S^{-1}F(\omega_{0})\rangle\mu(\lbrace\omega_{0}\rbrace)$ is not invertible in $\mathcal A$. Then $F:\Omega\setminus\lbrace\omega_{0}\rbrace \to U$ is not a continuous frame for Hilbert $C^{\ast}$-module $U$.
\end{theorem}

\begin{proof}
By the reconstruction formula we have
\begin{align*}
F(\omega_{0}) &=\int_{\Omega}\langle F(\omega_{0}),S^{-1}F(\omega_{0})\rangle F(\omega) d\mu(\omega)\\
&= \int_{\Omega\setminus\lbrace\omega_{0}\rbrace} \langle F(\omega_{0}) ,S^{-1}F(\omega)\rangle F(\omega) d\mu(\omega)+\langle F(\omega_{0}) ,S^{-1}F(\omega_{0})\rangle F(\omega_{0})\mu(\lbrace\omega_{0}\rbrace),
\end{align*}
then
\begin{equation*}
F(\omega_{0})(1_{\mathcal A}-\langle F(\omega_{0}) ,S^{-1}F(\omega_{0})\rangle\mu(\lbrace\omega_{0}\rbrace)) =\int_{\Omega\setminus\lbrace\omega_{0}\rbrace} \langle F(\omega_{0}) ,S^{-1}F(\omega)\rangle F(\omega) d\mu(\omega).
\end{equation*}
Now assume that $F:\Omega\setminus\lbrace\omega_{0}\rbrace \to U$ be a continuous frame for $U$. Since $S^{-\frac{1}{2}}$ is invertible, then $S^{-\frac{1}{2}}F:\Omega\setminus\lbrace\omega_{0}\rbrace \to U$ is also a continuous frame with the frame bounds $C,D>0$. Then for all $f\in U$,
\begin{equation*}
C\langle f,f\rangle \leq\int_{\Omega\setminus\lbrace\omega_{0}\rbrace} \langle f,S^{-\frac{1}{2}}F(\omega)\rangle\langle S^{-\frac{1}{2}}F(\omega),f\rangle d\mu(\omega)\leq D\langle f,f\rangle,
\end{equation*}
and for $f=S^{-\frac{1}{2}}F(\omega_{0})$ we have
\begin{align*}
C\langle S^{-\frac{1}{2}}F(\omega_{0}),S^{-\frac{1}{2}}F(\omega_{0})\rangle & \leq\int_{\Omega\setminus\lbrace\omega_{0}\rbrace} \langle S^{-\frac{1}{2}}F(\omega_{0}),S^{-\frac{1}{2}}F(\omega)\rangle\langle S^{-\frac{1}{2}}F(\omega),S^{-\frac{1}{2}}F(\omega_{0})\rangle d\mu(\omega)\\
& \leq D\langle S^{-\frac{1}{2}}F(\omega_{0}),S^{-\frac{1}{2}}F(\omega_{0})f\rangle,
\end{align*}
then
\begin{equation*}
C\langle F(\omega_{0}),S^{-1}F(\omega_{0})\rangle\leq\int_{\Omega\setminus\lbrace\omega_{0}\rbrace} \langle F(\omega_{0}),S^{-1}F(\omega)\rangle\langle S^{-1}F(\omega),F(\omega_{0})\rangle d\mu(\omega).
\end{equation*}
Define
\begin{align*}
\psi: & \Omega\longrightarrow\mathcal A\\
& \omega\longmapsto\langle F(\omega_{0}),S^{-1}F(\omega)\rangle
\end{align*}
By lemma \ref{psi}, $\psi\in L^{2}(\Omega ,\mathcal A)$ and 
\begin{equation*}
 \psi(\omega_{0})=\int_{\Omega\setminus\lbrace\omega_{0}\rbrace}\psi(\omega)\psi(\omega)^{*} d\mu(\omega)+(\psi(\omega_{0}))^{2}\mu(\lbrace\omega_{0}\rbrace).
\end{equation*}
Also
\begin{equation*}
C\psi(\omega_{0})\leq\int_{\Omega\setminus\lbrace\omega_{0}\rbrace}\psi(\omega)\psi(\omega)^{*} d\mu(\omega),
\end{equation*}
then
\begin{equation*}
C\psi(\omega_{0})\leq\psi(\omega_{0})-(\psi(\omega_{0}))^{2}\mu(\lbrace\omega_{0}\rbrace).
\end{equation*}
So $Ct\leq t-t^{2}\mu(\lbrace\omega_{0}\rbrace)$ holds for any $t$ in $\sigma(\psi(\omega_{0}))$, the spectrum of $\psi(\omega_{0})$.\\
Since $1_{\mathcal A}-\psi(\omega_{0})\mu(\lbrace\omega_{0}\rbrace)$ is not invertible and
\begin{equation*}
1_{\mathcal A}-\psi(\omega_{0})\mu(\lbrace\omega_{0}\rbrace)=(\frac{1}{\mu(\lbrace\omega_{0}\rbrace)}1_{\mathcal A}-\psi(\omega_{0}))\mu(\lbrace\omega_{0}\rbrace),
\end{equation*}
so $\frac{1}{\mu(\lbrace\omega_{0}\rbrace)}\in\sigma(\psi(\omega_{0}))$. Therefore,
\begin{equation*}
C\frac{1}{\mu(\lbrace\omega_{0}\rbrace)}\leq\frac{1}{\mu(\lbrace\omega_{0}\rbrace)} -\frac{1}{(\mu(\lbrace\omega_{0}\rbrace))^{2}}\mu(\lbrace\omega_{0}\rbrace)=0.
\end{equation*}
This is a contradiction and complete the proof.
\end{proof}     
%-------------------------------------------------------------------------------------------------------------------------------------------
We need the following lemma to prove the next theorems that  states an important property of the canonical dual frame of a given continuous frame for Hilbert $C^{\ast}$-modules.
%-------------------------------------------------------------------------------------------------------------------------------------------
\begin{lemma}\label{c-dual}
Let $F:\Omega \to U$ be a continuous frame for Hilbert $C^{\ast}$-module $U$ over a unital $C^*$-algebra $\mathcal A$ with frame operator $S$ and pre-frame operator $T$. Let $f\in U$ and suppose that $f=\int_{\Omega}\varphi(\omega) F(\omega) d\mu(\omega)$, for some $\varphi\in L^{2}(\Omega ,A)$. Then
\begin{align*}
\int_{\Omega}\vert\varphi(\omega)^{*}\vert^{2} d\mu(\omega)= & \int_{\Omega}\langle f,S^{-1}F(\omega)\rangle\langle S^{-1}F(\omega),f\rangle d\mu(\omega)\\
& +\int_{\Omega}(\varphi(\omega)-\langle f,S^{-1}F(\omega)\rangle)(\varphi(\omega)^{*}-\langle S^{-1}F(\omega),f\rangle) d\mu(\omega).
\end{align*}
\end{lemma}
\begin{proof}
For each $\omega\in\Omega$ we can write
\begin{equation*}
\varphi(\omega)= (\varphi(\omega)-\langle f,S^{-1}F(\omega)\rangle)+\langle f,S^{-1}F(\omega)\rangle
\end{equation*}
Since $F$ is a frame so, $f=\int_{\Omega}\langle f,S^{-1}F(\omega)\rangle F(\omega) d\mu(\omega)$. Then
\begin{align*}
& \int_{\Omega}(\varphi(\omega)-\langle f,S^{-1}F(\omega)\rangle) F(\omega) d\mu(\omega) =0\\
\Longrightarrow\; & \lbrace(\varphi(\omega)-\langle f,S^{-1}F(\omega)\rangle)\rbrace_{\omega\in\Omega}\in Ker(T)
\end{align*}
Also $\lbrace\langle f,S^{-1}F(\omega)\rangle\rbrace_{\omega\in\Omega}=\lbrace\langle S^{-1}F(\omega),f\rangle\rbrace_{\omega\in\Omega}\in R(T^{*})$  and $L^{2}(\Omega ,A)=Ker(T)\oplus R(T^{*})$.\\
It shows that
\begin{equation*}
\langle\lbrace\langle f,S^{-1}F(\omega)\rangle\rbrace_{\omega\in\Omega}, \lbrace(\varphi(\omega)-\langle f,S^{-1}F(\omega)\rangle)\rbrace_{\omega\in\Omega}\rangle =0.
\end{equation*}
Then
\begin{align*}
0 & =\int_{\Omega}\langle f,S^{-1}F(\omega)\rangle(\varphi(\omega)-\langle f,S^{-1}F(\omega)\rangle)^{*} d\mu(\omega)\\
& =\int_{\Omega}\langle f,S^{-1}F(\omega)\rangle\langle S^{-1}F(\omega),f\rangle d\mu(\omega) -\int_{\Omega}\langle f,S^{-1}F(\omega)\rangle\varphi(\omega)^{*} d\mu(\omega).
\end{align*}
Note that
\begin{align*}
\int_{\Omega}(\varphi(\omega)-\langle f,S^{-1}F(\omega)\rangle) & (\varphi(\omega)^{*}-\langle S^{-1}F(\omega),f\rangle) d\mu(\omega) =\\
& \int_{\Omega}\varphi(\omega)\varphi(\omega)^{*} d\mu(\omega) +\int_{\Omega}\varphi(\omega)\langle S^{-1}F(\omega),f\rangle d\mu(\omega)\\
& -\int_{\Omega}\langle f,S^{-1}F(\omega)\rangle\varphi(\omega)^{*} d\mu(\omega) +\int_{\Omega}\langle f,S^{-1}F(\omega)\rangle\langle S^{-1}F(\omega),f\rangle d\mu(\omega).
\end{align*}
Hence the proof is complete.
\end{proof}         
%-------------------------------------------------------------------------------------------------------------------------------------------
In the following, we generalize Theorem \ref{c-f-less-no} by removing a measurable subset of $\Omega$.  
%-------------------------------------------------------------------------------------------------------------------------------------------
\begin{theorem}
Let $F:\Omega \to U$ be a continuous frame for Hilbert $C^{\ast}$-module $U$ over a unital $C^*$-algebra $\mathcal A$ with frame operator $S$. Also suppose that $\Omega_{1}$ is a measurable subset of $\Omega$ such that $0<\mu(\Omega_{1})<\infty$ and $f=\int_{\Omega_{1}}F(\omega)d\mu(\omega)\neq 0$. If $\langle f,S^{-1}F(\omega)\rangle=\chi_{\Omega_{1}}(\omega)$ for all $\omega\in\Omega$, then the mapping $F:\Omega\setminus\Omega_{1} \to U$ is not a continuous frame for Hilbert $C^{\ast}$-module $U$.
\end{theorem}
\begin{proof}
We knowe,
\begin{equation*}
\int_{\Omega}\chi_{\Omega_{1}}(\omega)F(\omega)d\mu(\omega)=\int_{\Omega_{1}}F(\omega)d\mu(\omega)=f= \int_{\Omega}\langle f,S^{-1}F(\omega)\rangle(\omega)F(\omega)d\mu(\omega).
\end{equation*}
Also
\begin{equation*}
\int_{\Omega}(\chi_{\Omega_{1}}(\omega))(\chi_{\Omega_{1}}(\omega))^{*}d\mu(\omega)=\int_{\Omega}(\chi_{\Omega_{1}}(\omega))^{2}d\mu(\omega)=\mu(\Omega_{1}).
\end{equation*}
So by the Lemma \ref{c-dual} we have
\begin{align*}
\mu(\Omega_{1}) &=\int_{\Omega}\langle f,S^{-1}F(\omega)\rangle\langle S^{-1}F(\omega),f\rangle d\mu(\omega)\\
&+\int_{\Omega}(\chi_{\Omega_{1}}(\omega)-\langle f,S^{-1}F(\omega)\rangle)(\chi_{\Omega_{1}}(\omega)-\langle S^{-1}F(\omega),f\rangle)d\mu(\omega).
\end{align*}
Now if $\chi_{\Omega_{1}}(\omega)=\langle f,S^{-1}F(\omega)\rangle$ for all $\omega\in\Omega$, then
\begin{align*}
\mu(\Omega_{1}) &=\int_{\Omega_{1}}\langle f,S^{-1}F(\omega)\rangle\langle S^{-1}F(\omega),f\rangle d\mu(\omega)\\
&+\int_{\Omega\setminus\Omega_{1}}\langle f,S^{-1}F(\omega)\rangle\langle S^{-1}F(\omega),f\rangle d\mu(\omega)\\
&+\int_{\Omega}(\chi_{\Omega_{1}}(\omega)-\langle f,S^{-1}F(\omega)\rangle)(\chi_{\Omega_{1}}(\omega)-\langle S^{-1}F(\omega),f\rangle)d\mu(\omega)\\
&=\mu(\Omega_{1}) +\int_{\Omega\setminus\Omega_{1}}\langle f,S^{-1}F(\omega)\rangle\langle S^{-1}F(\omega),f\rangle d\mu(\omega) +0,
\end{align*}
then
\begin{equation*}
\langle\lbrace\langle f,S^{-1}F(\omega)\rangle\rbrace_{\omega\in\Omega\setminus\Omega_{1}},\lbrace\langle S^{-1}F(\omega),f\rangle\rbrace_{\omega\in\Omega\setminus\Omega_{1}}\rangle =\int_{\Omega\setminus\Omega_{1}}\langle f,S^{-1}F(\omega)\rangle\langle S^{-1}F(\omega),f\rangle d\mu(\omega)=0.
\end{equation*}
This shows that there exists a nonzero element $S^{-1}f\in U$ such that
\begin{equation*}
\langle S^{-1}f,F(\omega)\rangle =\langle F(\omega),S^{-1}f\rangle =0,
\end{equation*}
for all $\omega\in\Omega\setminus\Omega_{1}$. If $F:\Omega\setminus\Omega_{1} \to U$ is a continuous frame for $U$, then $Ker(T^{*}_{F\mid_{\Omega\setminus\Omega_{1}}})=\lbrace 0\rbrace$ and $f=0$, which is a contradiction.
\end{proof} 
%-------------------------------------------------------------------------------------------------------------------------------------------
The following theorem shows when does removing some elements from a continuous frame cause that the remaining set is a continuous frame.
%-------------------------------------------------------------------------------------------------------------------------------------------
\begin{theorem}\label{c-f-lees}
Let $F:\Omega \to U$ be a continuous frame for Hilbert $C^{\ast}$-module $U$ over a unital $C^*$-algebra $\mathcal A$ with frame operator $S$. Let $1_{\mathcal A}$ be the identity element of $\mathcal A$ and $\omega_{0}\in\Omega$ such that $1_{\mathcal A}-\langle F(\omega_{0}),S^{-1}F(\omega_{0})\rangle\mu(\lbrace\omega_{0}\rbrace)$ is invertible in $\mathcal A$. Then $F:\Omega\setminus\lbrace\omega_{0}\rbrace \to U$ is a continuous frame for Hilbert $C^{\ast}$-module $U$.
\end{theorem}
\begin{proof}
By reconstruction formula we have
\begin{align*}
F(\omega_{0}) &=\int_{\Omega}\langle F(\omega_{0}),S^{-1}F(\omega)\rangle F(\omega) d\mu(\omega)\\
& =\int_{\Omega\setminus\lbrace\omega_{0}\rbrace}\langle F(\omega_{0}),S^{-1}F(\omega)\rangle F(\omega) d\mu(\omega)+ \langle F(\omega_{0}),S^{-1}F(\omega_{0})\rangle F(\omega_{0})\mu(\lbrace\omega_{0}\rbrace)
\end{align*}
then
\begin{align*}
& (1_{\mathcal A}-\langle F(\omega_{0}),S^{-1}F(\omega_{0})\rangle\mu(\lbrace\omega_{0}\rbrace)) F(\omega_{0})= 
\int_{\Omega\setminus\lbrace\omega_{0}\rbrace}\langle F(\omega_{0}),S^{-1}F(\omega)\rangle F(\omega) d\mu(\omega)\\
\Longrightarrow\; & F(\omega_{0})=(1_{\mathcal A}-\langle F(\omega_{0}),S^{-1}F(\omega_{0})\rangle\mu(\lbrace\omega_{0}\rbrace))^{-1} \int_{\Omega\setminus\lbrace\omega_{0}\rbrace}\langle F(\omega_{0}),S^{-1}F(\omega)\rangle F(\omega) d\mu(\omega).
\end{align*}
Put $a:=(1_{\mathcal A}-\langle F(\omega_{0}),S^{-1}F(\omega_{0})\rangle\mu(\lbrace\omega_{0}\rbrace)) F(\omega_{0})^{-1}$, then by using Lemma \ref{NO}, we have
\begin{align*}
\langle f &,F(\omega_{0})\rangle\langle F(\omega_{0}),f\rangle\\
& =\langle f,a\int_{\Omega\setminus\lbrace\omega_{0}\rbrace}\langle F(\omega_{0}),S^{-1}F(\omega)\rangle F(\omega) d\mu(\omega)\rangle \langle a\int_{\Omega\setminus\lbrace\omega_{0}\rbrace}\langle F(\omega_{0}),S^{-1}F(\omega)\rangle F(\omega) d\mu(\omega) ,f\rangle\\
& =\langle f,\int_{\Omega\setminus\lbrace\omega_{0}\rbrace}\langle F(\omega_{0}),S^{-1}F(\omega)\rangle F(\omega) d\mu(\omega)\rangle a^{*}a\langle \int_{\Omega\setminus\lbrace\omega_{0}\rbrace}\langle F(\omega_{0}),S^{-1}F(\omega)\rangle F(\omega) d\mu(\omega) ,f\rangle\\
& \leq\Vert a\Vert^{2}\int_{\Omega\setminus\lbrace\omega_{0}\rbrace}\langle f,F(\omega)\rangle\langle S^{-1}F(\omega),F(\omega_{0})\rangle d\mu(\omega) \int_{\Omega\setminus\lbrace\omega_{0}\rbrace}\langle F(\omega_{0}),S^{-1}F(\omega),\rangle\langle F(\omega),f\rangle d\mu(\omega)\\
& =\Vert a\Vert^{2}\langle\lbrace\langle f,F(\omega)\rangle\rbrace_{\omega\in\Omega\setminus\lbrace\omega_{0}\rbrace} ,\lbrace\langle F(\omega_{0}),S^{-1}F(\omega)\rangle\rbrace_{\omega\in\Omega\setminus\lbrace\omega_{0}\rbrace}\rangle\\
&\qquad\qquad\qquad\qquad\qquad\langle\lbrace\langle F(\omega_{0}),S^{-1}F(\omega)\rangle\rbrace_{\omega\in\Omega\setminus\lbrace\omega_{0}\rbrace} ,\lbrace\langle f,F(\omega)\rangle\rbrace_{\omega\in\Omega\setminus\lbrace\omega_{0}\rbrace}\rangle\\
& \leq\Vert a\Vert^{2}\Vert \lbrace\langle F(\omega_{0}),S^{-1}F(\omega)\rangle\rbrace_{\omega\in\Omega\setminus\lbrace\omega_{0}\rbrace}\Vert^{2} \langle\lbrace\langle f,F(\omega)\rangle\rbrace_{\omega\in\Omega\setminus\lbrace\omega_{0}\rbrace} ,\lbrace\langle f,F(\omega)\rangle\rbrace_{\omega\in\Omega\setminus\lbrace\omega_{0}\rbrace}\rangle\\
& =\Vert a\Vert^{2}\Vert\int_{\Omega\setminus\lbrace\omega_{0}\rbrace}\langle F(\omega_{0}),S^{-1}F(\omega)\rangle\langle S^{-1}F(\omega),F(\omega_{0}) d\mu(\omega)\Vert\int_{\Omega\setminus\lbrace\omega_{0}\rbrace}\langle f,F(\omega)\rangle\langle F(\omega),f\rangle d\mu(\omega).
\end{align*}
Put
\begin{equation*}
 k:=\Vert a\Vert^{2}\Vert\int_{\Omega\setminus\lbrace\omega_{0}\rbrace}\langle F(\omega_{0}),S^{-1}F(\omega)\rangle\langle S^{-1}F(\omega),F(\omega_{0})\rangle d\mu(\omega)\Vert,
\end{equation*}
then
\begin{equation*}
 \langle f,F(\omega_{0})\rangle\langle F(\omega_{0}),f\rangle\leq k\int_{\Omega\setminus\lbrace\omega_{0}\rbrace}\langle f,F(\omega)\rangle\langle F(\omega),f\rangle d\mu(\omega).
\end{equation*}
Since $F$ is a continuous frame so there exsists constant $A>0$ such that
\begin{align*}
A\langle f,f\rangle & \leq\int_{\Omega}\langle f,F(\omega)\rangle\langle F(\omega),f\rangle d\mu(\omega)\\
& =\int_{\Omega\setminus\lbrace\omega_{0}\rbrace}\langle f,F(\omega)\rangle\langle F(\omega),f\rangle d\mu(\omega) +\langle f,F(\omega_{0})\rangle\langle F(\omega_{0}),f\rangle\mu(\lbrace\omega_{0}\rbrace)\\
& \leq\int_{\Omega\setminus\lbrace\omega_{0}\rbrace}\langle f,F(\omega)\rangle\langle F(\omega),f\rangle d\mu(\omega) +k\mu(\lbrace\omega_{0}\rbrace)\int_{\Omega\setminus\lbrace\omega_{0}\rbrace}\langle f,F(\omega)\rangle\langle F(\omega),f\rangle d\mu(\omega)\\
& =(1+k\mu(\lbrace\omega_{0}\rbrace))\int_{\Omega\setminus\lbrace\omega_{0}\rbrace}\langle f,F(\omega)\rangle\langle F(\omega),f\rangle d\mu(\omega).
\end{align*}
This implies that $F:\Omega\setminus\lbrace\omega_{0}\rbrace \to U$ satisfies the lower frame bound condition with bound $\frac{A}{1+k\mu(\lbrace\omega_{0}\rbrace)}$. 
Obviously $F:\Omega\setminus\lbrace\omega_{0}\rbrace \to U$ satisfies the upper frame bound condition. This completes the proof.
\end{proof}           
%-------------------------------------------------------------------------------------------------------------------------------------------
\begin{corollary}
Let $F:\Omega \to U$ be a continuous frame for Hilbert $C^{\ast}$-module $U$ over a unital $C^*$-algebra $\mathcal A$ with frame operator $S$. Let $1_{\mathcal A}$ be the identity element of $\mathcal A$ and $\omega_{0}\in\Omega$ such that $1_{\mathcal A}-\langle F(\omega_{0}),S^{-1}F(\omega_{0})\rangle\mu(\lbrace\omega_{0}\rbrace)$ is invertible in $\mathcal A$. Then $F$ is not a Riesz-type frame.
\end{corollary}
\begin{proof}
By the previous theorem, $F:\Omega\setminus\lbrace\omega_{0}\rbrace \to U$ is a continuous frame for $U$. Let $G:\Omega\setminus\lbrace\omega_{0}\rbrace \to U$ be the canonical dual frame for $F:\Omega\setminus\lbrace\omega_{0}\rbrace \to U$ and suppose that $G(\omega_{0})=0$. Then $S^{-1}F\neq G$ and for each $f,g\in U$,
\begin{align*}
\int_{\Omega}\langle f,G(\omega)\rangle\langle F(\omega),g\rangle d\mu(\omega) & =\int_{\Omega\setminus\lbrace\omega_{0}\rbrace}\langle f,G(\omega)\rangle\langle F(\omega),g\rangle d\mu(\omega) +\langle f,G(\omega_{0})\rangle\langle F(\omega_{0}),g\rangle\mu(\lbrace\omega_{0}\rbrace)\\
& = \langle f,g\rangle +0= \langle f,g\rangle.
\end{align*}
Hence $G:\Omega \to U$ is a dual frame for $F:\Omega \to U$ and $S^{-1}F\neq G$.
\end{proof}      
%-------------------------------------------------------------------------------------------------------------------------------------------
 By the same argument for the proof of the Theorem \ref{c-f-lees}, we can generalize it for an orbitrary dual $G$ of continuous frame $F$ as follows.
%-------------------------------------------------------------------------------------------------------------------------------------------
\begin{remark}
Let $F:\Omega \to U$ be a continuous frame for Hilbert $C^{\ast}$-module $U$ over a unital $C^*$-algebra $\mathcal A$ with the frame bounds $A_{F},B_{F}$. Also suppose that $1_{\mathcal A}$ be the identity element of $\mathcal A$ and $\omega_{0}\in\Omega$ such that $1_{\mathcal A}-\langle F(\omega_{0}),G(\omega_{0})\rangle\mu(\lbrace\omega_{0}\rbrace)$ is invertible in $\mathcal A$ where $G:\Omega \to U$ is an orbitrary dual of $F$ with the frame bounds $A_{G},B_{G}$. Then $F:\Omega\setminus\lbrace\omega_{0}\rbrace \to U$ is a continuous frame for Hilbert $C^{\ast}$-module $U$ with the lower frame bound $\frac{A_{F}}{1+\Vert a\Vert^{2}B_{G}\Vert F(\omega_{0})\Vert^{2}\mu(\lbrace\omega_{0}\rbrace)}$.
\end{remark}
%-------------------------------------------------------------------------------------------------------------------------------------------


\begin{thebibliography}{99}
\bibitem{CH} O. Christensen, \emph{An introduction to frames and Riesz bases}, Birkhauser, Boston, 2016.

\bibitem{DG} I. Daubechies, A. Grassman and Y. Meyer, \emph{Painless nonothogonal expanisions}, J. Math. Phys.,  
\textbf{27} (1986), 1271--1283.

\bibitem{DS} R.J. Duffin and A.C. Schaeffer, \emph{A class of nonharmonic Fourier series}, Trans. Amer.
Math. Soc., \textbf{72} (1952), 341--366.

\bibitem{FL} M. Frank and D. R. Larson, \emph{Frames in Hilbert $C^*$-modules and $C^*$-algebras}, J. Operator Theory, \textbf{48} (2002), 273--314.

\bibitem{GHSH1} H. Ghasemi and T.L. Shateri, \emph{Continuous $*$-controlled frames in Hilbert $C^*$-modules}, Caspian J. Math. Scien., \textbf{11}(2) (2022), 448--460, DOI: 10.22080/cjms.2022.21850.1590.

\bibitem{GHSH} H. Ghasemi and T.L. Shateri, \emph{Continuous Frames in Hilbert $C^*$-Modules}, arXiv preprint arXiv:2208.06799 (2022).

\bibitem{JI} W. Jing, \emph{Frames in Hilbert $C^*$-modules}, Ph.D. Thesis, University of Central Florida Orlando, Florida, 2006.

\bibitem{LAN} E.C. Lance, \emph{Hilbert $C^*$--Modules: A Toolkit for Operator Algebraist}, 144 pages, vol. 210 of London Mathematical Society Lecture Note Series, Cambridge University Press, Cambridge, UK, (1995).

\bibitem{MA}  B. Magajna, \emph{Hilbert $C^*$-modules in which all closed submodules are complemented}, Proc. Amer.
Math. Soc., \textbf{125} (3) (1997), 849--852.

\bibitem{MT} V. M. Manuilov and E. V. Troitsky, \emph{Hilbert $C^*$-Modules: Translations of mathematical monographs}, American Mathematical Soc, (2005).

\bibitem{RR} M. Rashidi-Kouchi and A. Rahimi, \emph{On controlled frames in Hilbert $C^*$-modules}, Int. J. Wavelets Multiresolut.
Inf. Process., \textbf{15}(4) (2017), 1750038 (15 pages), DOI: 10.1142/S0219691317500382.

\bibitem{SH} T.L. Shateri, \emph{$*$-controlled frames in Hilbert $C^*$-modules}, Inter. J. Wavelets Multiresolut. Inf. Process., \textbf{19}(03) 2050080 (2021), DOI: 10.1142/S0219691320500800.

\bibitem{XZ} X-C Xiao, X-M Zeng. \emph{Some properties of g-frames in Hilbert $C^*$-modules}, Journal of Mathematical Analysis and Applications \textbf{363}(2) (2010), 399--408.

\bibitem{OLS} N.E. Wegge-Olsen, \emph{K-theory and $C^*$-algebras: a friendly approach}, Oxford University Press, (1993).
\end{thebibliography}
\end{document}